\documentclass{article}
\usepackage[utf8]{inputenc}
\usepackage[T1]{fontenc}
\usepackage{tabularx}
\usepackage{array}
\usepackage{amsmath,amsfonts,stmaryrd,amssymb, amsthm} % Math packages

\usepackage{enumerate} % Custom item numbers for enumerations
\usepackage{enumitem}
\usepackage{multicol}
\usepackage{booktabs}

\usepackage[framemethod=tikz]{mdframed} % Allows defining custom boxed/framed environments

\usepackage{listings} % File listings, with syntax highlighting
\usepackage{hyperref}
\hypersetup{
	colorlinks=true,
	linkcolor=blue,
	filecolor=magenta,      
	urlcolor=cyan,
}

\usepackage{xfrac}%required for big quotient symbol \sfrac

\usepackage[title]{appendix}%required for using appendix

\usepackage{graphicx} % Required for including images
\graphicspath{{./images/}} % Set the default folder for images
\usepackage{caption}
\usepackage{subcaption}

\theoremstyle{plain}

\theoremstyle{definition}
\numberwithin{equation}{section}

\usepackage{algorithm}
\usepackage{algpseudocode}

\usepackage{geometry} % Required for adjusting page dimensions and margins

\usepackage[utf8]{inputenc} % Required for inputting international characters
\usepackage[T1]{fontenc} % Output font encoding for international characters

\usepackage{XCharter} % Use the XCharter fonts

\mdfdefinestyle{commandline}{
	leftmargin=10pt,
	rightmargin=10pt,
	innerleftmargin=15pt,
	middlelinecolor=black!50!white,
	middlelinewidth=2pt,
	frametitlerule=false,
	backgroundcolor=black!5!white,
	frametitle={Command Line},
	frametitlefont={\normalfont\sffamily\color{white}\hspace{-1em}},
	frametitlebackgroundcolor=black!50!white,
	nobreak,
}

\mdfdefinestyle{file}{
	innertopmargin=1.6\baselineskip,
	innerbottommargin=0.8\baselineskip,
	topline=false, bottomline=false,
	leftline=false, rightline=false,
	leftmargin=2cm,
	rightmargin=2cm,
	singleextra={%
		\draw[fill=black!10!white](P)++(0,-1.2em)rectangle(P-|O);
		\node[anchor=north west]
		at(P-|O){\ttfamily\mdfilename};
		\def\l{3em}
		\draw(O-|P)++(-\l,0)--++(\l,\l)--(P)--(P-|O)--(O)--cycle;
		\draw(O-|P)++(-\l,0)--++(0,\l)--++(\l,0);
	},
	nobreak,
}

\mdfdefinestyle{question}{
	innertopmargin=1.2\baselineskip,
	innerbottommargin=0.8\baselineskip,
	roundcorner=5pt,
	nobreak,
	singleextra={%
		\draw(P-|O)node[xshift=1em,anchor=west,fill=white,draw,rounded corners=5pt]{%
		Question \theQuestion\questionTitle};
	},
}

\newcounter{Question} % Stores the current question number that gets iterated with each new question

\mdfdefinestyle{warning}{
	topline=false, bottomline=false,
	leftline=false, rightline=false,
	nobreak,
	singleextra={%
		\draw(P-|O)++(-0.5em,0)node(tmp1){};
		\draw(P-|O)++(0.5em,0)node(tmp2){};
		\fill[black,rotate around={45:(P-|O)}](tmp1)rectangle(tmp2);
		\node at(P-|O){\color{white}\scriptsize\bf !};
		\draw[very thick](P-|O)++(0,-1em)--(O);%--(O-|P);
	}
}

\mdfdefinestyle{info}{%
	topline=false, bottomline=false,
	leftline=false, rightline=false,
	nobreak,
	singleextra={%
		\fill[black](P-|O)circle[radius=0.4em];
		\node at(P-|O){\color{white}\scriptsize\bf i};
		\draw[very thick](P-|O)++(0,-0.8em)--(O);%--(O-|P);
	}
}

\usepackage{authblk}
\title{Graph Neural Networks and 3-Dimensional Topology}
\author[1]{Pavel Putrov}
\author[1,2]{Song Jin Ri}
\affil[1]{ICTP, Strada Costiera 11, Trieste 34151, Italy}
\affil[2]{SISSA, Via Bonomea 265, Trieste 34136, Italy}
\date{}                     %% if you don't need date to appear
\begin{document}
	
\renewcommand\thefootnote{\textcolor{red}{\arabic{footnote}}}

\def\Z{{\mathbb{Z}}}
\def\C{{\mathbb{C}}}
\def\R{{\mathbb{R}}}

\maketitle % Print the title

%----------------------------------------------------------------------------------------
%	INTRODUCTION
%----------------------------------------------------------------------------------------

\begin{abstract}
We test the efficiency of applying Geometric Deep Learning to the problems in low-dimensional topology in a certain simple setting. Specifically, we consider the class of 3-manifolds described by plumbing graphs and use Graph Neural Networks (GNN) for the problem of deciding whether a pair of graphs give homeomorphic 3-manifolds. We use supervised learning to train a GNN that provides the answer to such a question with high accuracy. Moreover, we consider reinforcement learning by a GNN to find a sequence of Neumann moves that relates the pair of graphs if the answer is positive.  The setting can be understood as a toy model of the problem of deciding whether a pair of Kirby diagrams give diffeomorphic 3- or 4-manifolds.
\end{abstract}

%\tableofcontents

\section{Introduction and Summary}
\label{section_introduction}

Geometric Deep Learning (GDL) \cite{bronstein2021geometric} is an area of Machine Learning (ML) that has been under very active development during the last few years. It combines various approaches to ML problems involving data that has some underlying geometric structure. The neural networks used in GDL are designed to naturally take into account the symmetries and the locality of the data. It has been successfully applied to problems involving computer vision, molecule properties, social or citation networks, particle physics, etc (see \cite{cao2020comprehensive} for a survey). It is natural to apply GDL techniques also to mathematical problems in topology. In general, ML has been already used in various problems in low-dimensional topology, knot theory in particular, \cite{hughes2020neural,Jejjala:2019kio,Gukov:2020qaj,davies2021advancing,kauffman2022rectangular,Craven:2021ckk,vernitski2022reinforcement,khan2021untangling,lisitsa2023supervised,gukov2023searching}, as well as various physics-related problems in geometry (for a recent survey see \cite{He:2023csq}).
However, the used neural network models were mostly not specific to GDL. 

The goal of this paper is to test the efficiency of GDL in a very simple setting in low-dimensional topology. Namely, we consider a special class of 3-manifolds known as plumbed, or graph, 3-manifolds. Those are 3-manifolds that are specified by a choice of a graph with particular features assigned to edges and vertices. Such 3-manifolds are therefore very well suited for analysis by Graph Neural Networks (GNN). GNN is one of the most important and used types of neural networks used in GDL. In general, GNN are designed to process data represented by graphs. 

In this paper, we use GNNs for the following problems involving plumbed 3-manifolds. Different (meaning not isomorphic) graphs can correspond to equivalent, i.e. homeomorphic, 3-manifolds. Note that in 3 dimensions (or less) any topological manifold has a unique smooth structure and the notions of homeomorphism and diffeomorphism are equivalent. It is known that a pair of graphs that produce two equivalent 3-manifolds must be related by a sequence of certain \textit{moves}, commonly known as Neumann moves \cite{neumann1981calculus}. These moves establish a certain equivalence relation on the graphs (in addition to the standard graph isomorphism). First, we consider a neural network that, as the input has a pair of plumbing graphs, and, as the output gives the decision whether the graphs correspond to homeomorphic 3-manifolds or not, i.e. whether the two graphs are equivalent, or not, in the sense described above. Supervised Learning (SL) is then used to train the network. The training dataset consists of randomly generated graph pairs, for which it is known whether the corresponding 3-manifolds are homeomorphic or not. The trained neural network, up until the very last layer, can be understood to produce an approximate topological invariant of plumbed 3-manifolds. 

Second, we consider a neural network for which the input is a plumbing graph and the output is a sequence of Neumann moves that ``simplifies'' the graph according to a certain criterion. The aim is to build a neural network such that if it is applied to equivalent graphs it simplifies them to the same graph. If the result is successful this can be used to provide an explicit demonstration that a given pair of graphs give two homeomorphic 3-manifolds. Reinforcement Learning (RL) is used to train the network. For both cases, SL and RL, we consider different architectures of the neural networks and compare their performance. 

Note that in principle there is an algorithm for determining whether two plumbing graphs give homeomorphic 3-manifolds or not, which was already presented in \cite{neumann1981calculus}. It involves bringing both graphs to a certain normal form (which is, in a sense, similar to the ``simplification'' process in the RL setup mentioned above) and then checking that normal forms are the same (i.e. isomorphic graphs). However, it is known that just checking isomorphism of graphs already goes beyond polynomial time. The plumbing graphs can be considered as a particular class of more general Kirby diagrams that can be used to describe arbitrary closed oriented 3-manifolds, with Neumann moves being generalized to the so-called Kirby moves. Even in this case, in principle there exists an algorithm of checking whether two Kirby diagrams produce homeomorphic 3-manifolds or not \cite{kuperberg2019algorithmic}. There is a also of version of Kirby diagrams and moves for smooth 4-manifolds. Moreover, in this case, however, an algorithm for the recognition of diffeomorphic pairs does not exist. In 4 dimensions the notion of diffeomorphism and homeomorphism are not the same. In particular, there exist pairs of manifolds that are homeomorphic but not diffeomorphic. While the classification of 4-manifolds up to homeomorphisms (with certain assumptions on the fundamental group) is relatively not difficult, classification up to diffeomorphisms is an important open question. The setup with plumbed 3-manifolds that we consider in this paper can be understood as a toy model for the problem of recognition of diffeomorphic pairs of general 3- and 4-manifolds, for which one can try to apply neural networks with similar architecture in the future.

The rest of the paper is organized as follows. In Section \ref{section_preliminaries} we review basic preliminaries about plumbed 3-manifolds and Graph Neural Networks needed for the analysis that follows. In Section \ref{section_supervised} we consider various GNN architectures for supervised learning of whether a pair of plumbing graphs provide homeomorphic 3-manifolds or not. In Section \ref{section_reinforcement} we consider reinforcement learning of the process of simplification of a plumbing graph representing a fixed (up to a homeomorphism) 3-manifold. Finally, we conclude with Section \ref{section_conclusion} where we discuss the obtained results and mention possible further directions. The Appendix \ref{appendix_algorithms} contains some basic algorithms that are specific to the problems considered in this paper. 
\section{Preliminaries}
\label{section_preliminaries}

\subsection{Plumbed 3-manifolds}

\begin{figure}[t]
	\centering
	\includegraphics[scale=2]{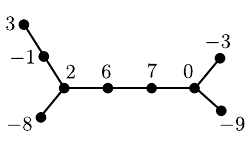}
	\caption{An example of a plumbing graph.}
	\label{fig_plumbing_example}
\end{figure}

In this section we review basic facts about plumbed 3-manifolds, also known as graph 3-manifolds. For a more detailed exposition we refer to the original paper \cite{neumann1981calculus}. First, let us describe how to build a 3-manifold from a \textit{plumbing graph}, or simply a \textit{plumbing}. For convenience, we restrict ourselves to the case when the graph is a tree, i.e., the graph is connected and acyclic. We will also consider the case of genus zero plumbings only. In this setting, apart from the graph itself, the only additional information that one needs to specify is the set of integer \textit{weights} $w(v)\in \mathbb{Z}$ labeling vertices $v\in V$ ($V$ denotes set of all vertices of the graph), also referred to as \textit{framings} in the context of topology. A typical plumbing graph looks like the one shown in Figure \ref{fig_plumbing_example}. The weights $w(v)$, together with standard graph data can be naturally encoded in an $|V|\times |V|$ matrix $a$ with integral elements $a_{ij}$ as follows. Outside of the diagonal this matrix coincides with the standard adjacency matrix of the graph (i.e. $a_{ij}=1$ if $i\neq j\in V$ are connected by an edge, and $a_{ij}=0$ otherwise). The diagonal elements are given by the weights: $a_{ii}=w(i)$.

One can build a 3-manifold corresponding to such a plumbing graph as follows. First, consider a graph containing a single vertex with weight $p\in \Z$ and no edges. To such a one-vertex graph we associate lens space 3-manifold $L(|p|,\pm 1)$, where the sign of $\pm 1$ coincides with the sign of $p$. It can be described as a quotient of the standard unit 3-sphere $S^3=\{|z_1|^2+|z_2|^2=1\,|\,(z_1,z_2)\in \C^2\}\subset \R^4\cong \C^2$ with respect to the action of cyclic group $\Z_{|p|}$ or order $|p|$, generated by $(z_1,z_2)\rightarrow (z_1e^{\frac{2\pi i}{|p|}},z_2e^{\mp\frac{2\pi i}{|p|}})$ transformation. This 3-manifold can be equivalently understood as a circle fibration over $S^2$ base with Euler number $p$. More explicitly, it can be constructed as follows. Let us start with two copies of $D^2\times S^1$ (where $D^2$ denotes 2-dimensional disk), that can be viewed as trivial circle fibrations over $D^2$. We then can glue two $D^2$'s along the common boundary $\partial D^2\cong S^1$ into $S^2$ (so that each $D^2$ can be understood as a hemisphere), with the $S^1$ fibers along the two boundaries being glued with relative rotation specified by a certain map $f:\partial D^2\cong S^1 \rightarrow SO(2)\cong S^1$. The homotopy class of such a map is completely determined by the ``winding number''. The homeomorphism class of the resulting closed 3-manifold only depends on this number. To obtain the lens space $L(p,1)$ one takes the winding number to be $p$. 

Next, consider a vertex with weight $p$ being a part of a general tree plumbing (as, for example, the one shown in Figure \ref{fig_plumbing_example}). For each edge coming out of the vertex, we remove a single $S^1$ fiber (over some generic point in the $S^2$ base) of the fibration together with its tubular neighborhood. The neighborhood can be chosen to be the restriction of the fibration to a small disk in the $S^2$ base that contains the chosen point. Such an operation, out of the original lens space $L(p,1)$, produces  a 3-manifold that has a boundary component $\partial D^2 \times S^1\cong S^1\times S^1=T^2$ for each edge coming out of the vertex. Having an edge between a pair of vertices in the graph then corresponds to gluing two $T^2$ boundary components in the way that the two circles, the fiber $S^1$ and the boundary $S^1$ of the small disk on the base, are swapped (with the orientation of one of the circles reversed, so that the resulting 3-manifold is orientable). Performing such operations to all the vertices and edges of the graph one obtains a 3-manifold that has no boundary components. This is the 3-manifold that one associates to the plumbing graph.

Equivalently, to a plumbing graph one can associate Dehn surgery diagram, in the way that each vertex $v\in V$ corresponds to an unknot framed by $w(v)$ and the presence of an edge between two vertices signifies that the corresponding unknots form a Hopf link.

\begin{figure}[t]
	\centering
    \begin{subfigure}[t]{0.31\textwidth}
        \centering
        \includegraphics[scale=1.7]{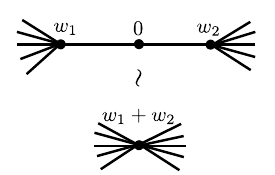}
        \caption{A Neumann move of type (a)}
    \end{subfigure}
    \hfill
	\begin{subfigure}[t]{0.31\textwidth}
	    \centering
        \includegraphics[scale=1.7]{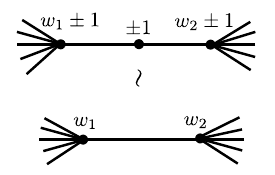}
        \caption{Neumann moves of type (b)}
	\end{subfigure}
    \hfill
    \begin{subfigure}[t]{0.31\textwidth}
	    \centering
        \includegraphics[scale=1.7]{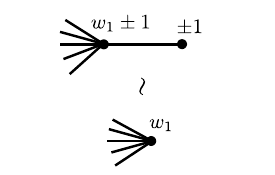}
        \caption{Neumann moves of type (c)}
	\end{subfigure}
    \hfill    
	\caption{There are 3 different types of Neumann moves, which preserve the resulting 3-manifold up to homeomorphism. Considering the sign, one can count 8 Neumann moves. Among them, there are 5 blow-up moves (by which a new vertex is created) and 3 blow-down moves (by which one or two vertices are annihilated).}
	\label{fig_nmoves}
\end{figure}

Applying the prescription described above to different graphs may result in homeomorphic 3-manifolds. In \cite{neumann1981calculus} it was proved that this happens if and only if the graphs can be related by a sequence of local graph transformations, or \textit{moves}, now commonly known as \textit{Neumann moves}, shown in Figure \ref{fig_nmoves}.

\subsection{Graph neural networks}
\label{subsection_preliminary_gnns}
%Graphs have a special data structure which models a set of finite or infinite objects (nodes) and pairwise relations (edges) between objects. 
Here we provide a brief review on some of GNNs for a later purpose. There are 3 main computational modules to build a typical GNN architecture: propagation modules, sampling modules and pooling modules. Since in this paper we will use only convolution operators, which are one of the most frequently used propagation modules, we focus on some of convolution operators. For a broad review on various modules, we refer the reader to \cite{ZHOU202057}.

Convolution operators are motivated by convolutional neural networks (CNN), which have achieved a notable progress in various areas. In general, the role of convolution operators can be described as
\begin{equation*}
    \mathbf{x}_i^{(k)}=\gamma^{(k)}\left(\mathbf{x}_i^{(k-1)}, \bigoplus_{j\in\mathcal{N}(i)}\phi^{(k)}\left(\mathbf{x}_i^{(k-1)}, \mathbf{x}_j^{(k-1)}, \mathbf{e}_{j,i}\right)\right),
\end{equation*}
where $\mathbf{x}_i^{(k)}\in\mathbb{R}^{F}$ denotes node features of node $i$ in the $k$-th layer and $\mathbf{e}_{j,i}\in\mathbb{R}^D$ denotes edge features of the edge connecting from node $j$ to node $i$. We also note that $\bigoplus$ over a neighborhood $\mathcal{N}(i)$ of node $i$ is a differentiable, permutation invariant function such as sum, mean and max, and $\gamma$ and $\phi$ denote differentiable functions such as Multi Layer Perceptrons (MLPs).

Among various convolution operators existing in the literature, the following will appear in the next sections.
\begin{itemize}
    \item Graph Embedding Network (GEN) \cite{li2019graph}\\ 
    GEN is designed for deep graph similarity learning and embeds each graph into a vector, called a graph embedding. More explicitly, it first computes initial node embeddings $\mathbf{x}^{(1)}_i$ from the node features $\mathbf{x}^{(0)}_i$ through MLP
    \[
        \mathbf{x}^{(1)}_i = \text{MLP}\left(\mathbf{x}^{(0)}_i\right),
    \]
    then it executes the single message propagation to compute node embeddings $\mathbf{x}^{(2)}_i$ by the information in its local neighbourhood $\mathcal{N}(i)$ \footnote{It is also possible to apply a finite number of propagation process iteratively, but we will only consider single propagation here.}
    \[
        \mathbf{x}^{(2)}_i = \text{MLP}\left(\mathbf{x}^{(1)}_{i}, \sum_{j\in\mathcal{N}(i)}\text{MLP}\left(\mathbf{x}^{(1)}_{i}, \mathbf{x}^{(1)}_{j}\right)\right).
    \]
    Once the node embeddings $ \mathbf{x}^{(2)}_i $ are computed, an aggregator computes a graph embedding by aggregating the set of node embeddings. In Section \ref{subsection_supervised_models}, we describe the details of the aggregator which we will apply not only to GEN but also the other models GCN and GAT.
    \item Graph Convolutional Network (GCN)\\
    GCN is introduced in \cite{kipf2017semisupervised} as a variant of convolutional neural networks for graphs. It operates as the following formula:
    \[
        \mathbf{z}_i=\mathbf{\Theta}^\intercal\sum_{j\in \mathcal{N}(i)\cup \{i\}}\frac{1}{\sqrt{\hat{d}_j\hat{d}_i}}\mathbf{x}_j,
    \]
    where $\mathbf{z}_i$ is the output for the $i$-th node, $\mathbf{\Theta}$ is a matrix of filter parameters, and $\hat{d}_i$ is the degree of $i$-th node.
    \item Graph Attention Network (GAT)\\
    GAT is proposed in \cite{velickovic2018graph}, which incorporates the attention mechanism into the message propagation. The mechanism of GAT can be formulated as
    \[
        \mathbf{x}'_i = \alpha_{i,i}\mathbf{\Theta}\mathbf{x}_i+\sum_{j\in\mathcal{N}(i)}\alpha_{i,j}\mathbf{\Theta}\mathbf{x}_j.
    \]
    Here the attention coefficients $\alpha$ are given by
    \[
        \alpha_{i,j}=\frac{\exp\left(\text{LeakyReLU}(\mathbf{a}^\intercal[\mathbf{\Theta}\mathbf{x}_i\Vert \mathbf{\Theta}\mathbf{x}_j)]\right)}{\sum_{k\in\mathcal{N}(i)\cup\{i\}}\exp\left(\text{LeakyReLU}(\mathbf{a}^\intercal[\mathbf{\Theta}\mathbf{x}_i\Vert \mathbf{\Theta}\mathbf{x}_k)]\right)},
    \]
    where the attention mechanism $\mathbf{a}$ is implemented by a single-layer feedfoward neural network, and $\Vert$ is the concatenation operator.
    \end{itemize}

    All the neural networks including GNNs are implemented based on PyTorch \cite{paszke2017pytorch} and PyTorch Geometric \cite{FeyLenssen2019pytorchgeometric}. \footnote{Python code is available on \href{https://github.com/songjin91/LearningPlumbings/tree/main}{Github}.}
\section{Supervised Learning}
\label{section_supervised}
In this section we use supervised learning to decide whether or not two plumbing graphs represent a same plumbed 3-manifold. We build 3 models GEN+GAT, GCN+GCN and GCN+GAT and examine their performance for the task.\footnote{Since our task is a graph similarity learning, we had expected that GEN could perform well and we performed GEN combined with various GNNs known in the literature. GEN+GAT performed best among different architectures we tried, and two models GCN+GCN and GCN+GAT are used for a benchmark because GCN and GAT are most commonly used GNNs.}

\subsection{Models}
\label{subsection_supervised_models}
All the models are designed to have two convolution operators, one aggregation layer and one classification layer. The models are named by concatenating the names of two convolution operators. For a fair comparison, we use the common aggregation layer and the classification layer, and all the layers have the same dimensions for both input and output. 

Since we have already reviewed the convolution operators in Section \ref{subsection_preliminary_gnns}, let us now elaborate on the common aggregation layer and classification layer. The aggregator computes a graph embedding by aggregating all of its node embeddings, passed from convolution operators. We use the aggregation layer proposed in \cite{li2017gated}, which is formulated by
\[
    \mathbf{h}_G = \text{MLP}_G\left(\sum_{i\in V}\text{Softmax}(\text{MLP}_{\text{gate}}(\mathbf{x}_i))\odot\text{MLP}(\mathbf{x}_i)\right),
\]
where $\mathbf{h}_G$ is a graph-level output and $\odot$ denotes element-wise multiplication. 

The classification layer plays a role to determine, for a given pair of plumbing graphs, whether or not they are equivalent. This layer has the concatenation of two graph embeddings as its input and classifies into two classes, class 0 and class 1. Here class 1 means two plumbing graphs are equivalent while class 0 denotes they are inequivalent. We implement the classification layer by using MLP with two hidden layers.

A detailed information of the architecture for 3 models are presented in Table \ref{table_sl_models}. For each layer in the table, the first element in the bracket followed by a name of model denotes the dimension of input vectors of the layer while the second one denotes the dimension of output embedding. 

\begin{table}[!htbp]
    \centering
    \caption{The architecture of 3 models with parameter values.}
    \label{table_sl_models}
    {\renewcommand{\arraystretch}{1.5}
    \begin{tabular}{|>{\centering\arraybackslash} m{0.2\textwidth}|>{\centering\arraybackslash} m{0.2\textwidth}|>{\centering\arraybackslash} m{0.2\textwidth}|>{\centering\arraybackslash} m{0.2\textwidth}|}        
        \hline
        Layers & GEN+GAT & GCN+GAT & GCN+GCN\\ 
        \hline
        First convolution & GEN(1, 128) &  \multicolumn{2}{c|}{GCN(1, 128)}\\
        \hline
        Second convolution & \multicolumn{2}{c|}{GAT(128, 128)} & GCN(128, 128)\\
        \hline
        Aggregation& \multicolumn{3}{c|}{Aggregator(128, 32)}\\
        \hline
        Classification & \multicolumn{3}{c|}{MLP(64, 2)}\\
        \hline
    \end{tabular}}
\end{table}

\subsection{Experimental Settings}
For training and validation, we put together datasets including 80,000 random pairs of plumbings generated by algorithms presented in Appendix \ref{appendix_algorithms}. More explicitly, the datasets consists of
\begin{itemize}
    \item 40,000 pairs of equivalent plumbings generated by \textsc{EquivPair}, Algorithm \ref{alg_equivpair}, with $N_{\max}=40$.\\ 
    To generate a pair of equivalent plumbings, the algorithm starts with a random plumbing created by \textsc{RandomPlumbing}, Algorithm \ref{alg_randomplumbing}, and iteratively applies Neumann moves using \textsc{RandomNeumannMove}, Algorithm \ref{alg_RandomNeumannMove}, up to $N_{\max}$ times, to each plumbing in the pair.
    \item 30,000 pairs of inequivalent plumbings generated by \textsc{InequivPair}, Algorithm \ref{alg_inequivpair}, with $N_{\max}=40$.\\
    It has a similar process to \textsc{EquivPair}, but it starts with a pair of inequivalent plumbings, each of which is separately generated by \textsc{RandomPlumbing}.\footnote{We note that two plumbings, generated by running \textsc{RandomPlumbing} twice, could be accidentally equivalent and this might affect the accuracy of models in training. However, we will ignore this since it is statistically insignificant.}
    \item 10,000 pairs of inequivalent plumbings generated by \textsc{TweakPair}, Algorithm \ref{alg_tweakpair}, with $N_{\max}=40$.\\ 
    This algorithm generates a pair of inequivalent plumbings, one of which is obtained by tweaking the other. Here, by tweaking a plumbing, we mean that we make a small change of the weight (or node feature) of a randomly chosen node in the plumbing. Since tweaking is different from Neumann moves, this process creates an inequivalent plumbing to the original one. After tweaking, it also applies \textsc{RandomNeumannMove} iteratively up to $N_{\max}$. These pairs are added into the datasets in order for the models to make the decision boundary more accurate, since for a pair generated by \textsc{InequivPair}, two plumbings might be quiet different due to random generators.
\end{itemize}

We divide the datasets into training and validation sets by the ratio 8:2. We train our models on training sets containing 64,000 pairs of plumbings up to 150 epochs. For each model, we use cross-entropy loss for a loss function and  Adam for an optimizer with the learning rate 0.001.

\begin{figure}
     \centering
     \begin{subfigure}[t]{0.48\textwidth}
         \centering
         \includegraphics[width=\textwidth]{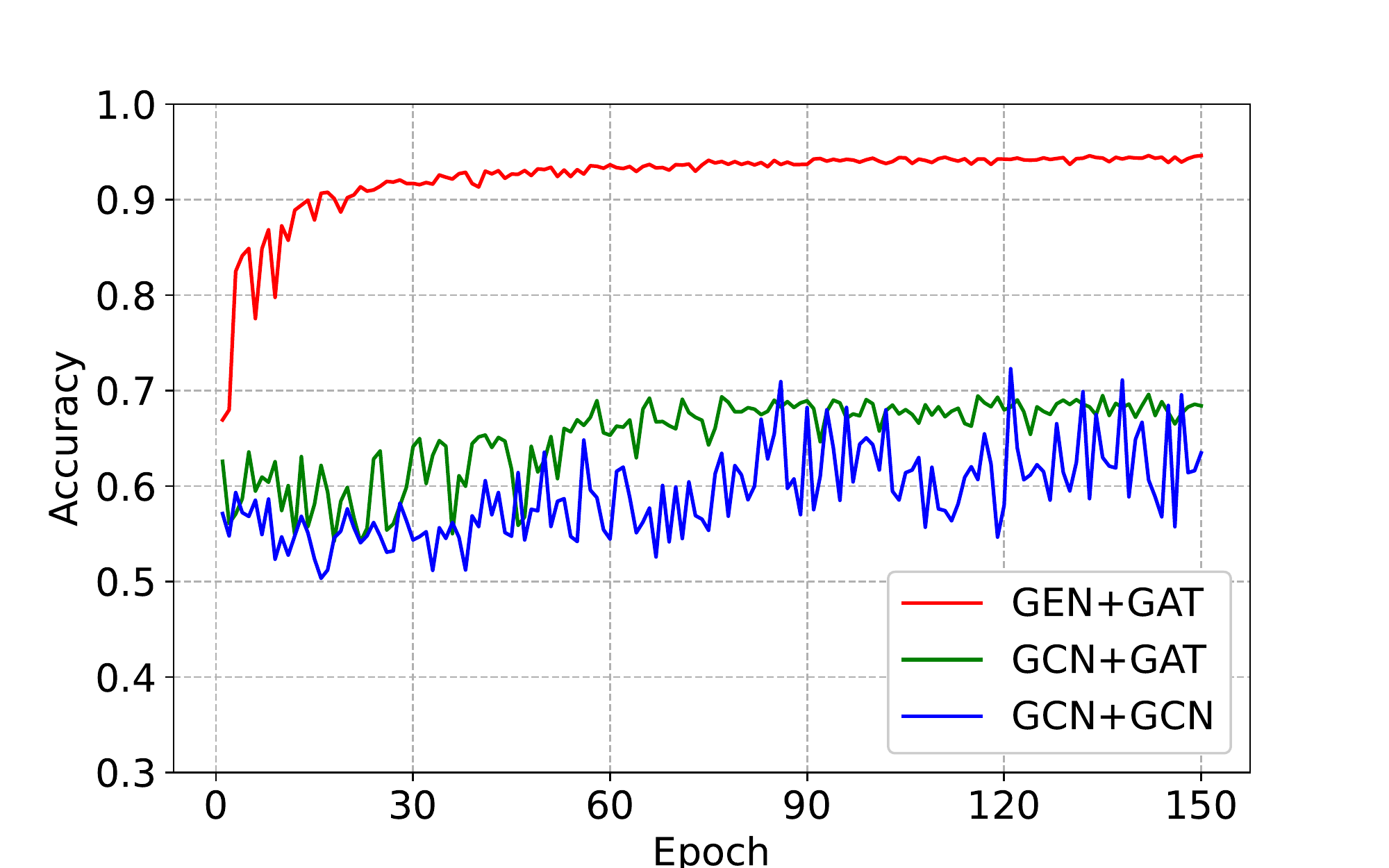}
         \caption{Performance comparison on validation sets.}
         \label{fig_performance_training}
     \end{subfigure}
     \hfill
     \begin{subfigure}[t]{0.48\textwidth}
         \centering
         \includegraphics[width=\textwidth]{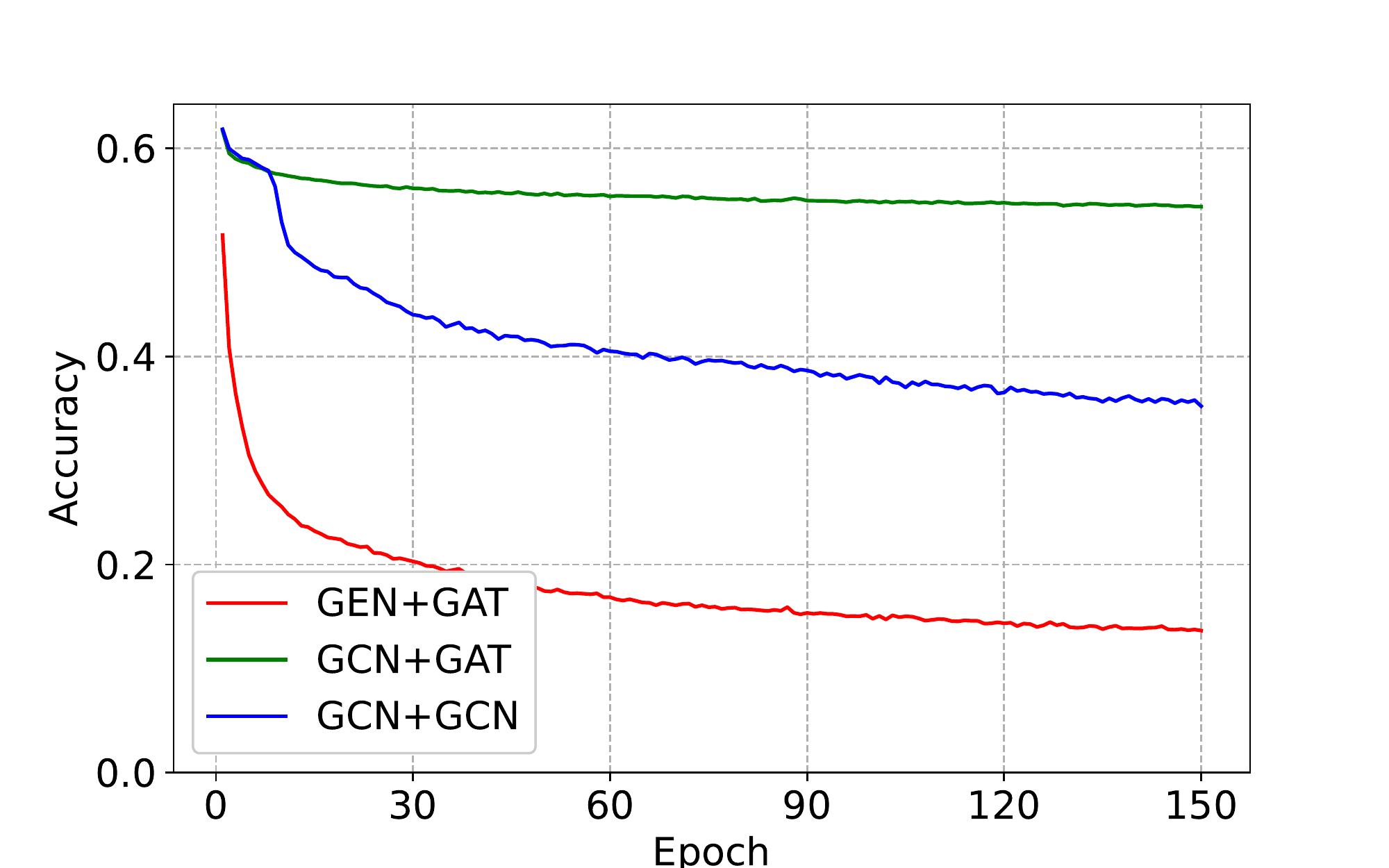}
         \caption{Loss comparison on validation sets.}
         \label{fig_loss_comparison}
    \end{subfigure}
    \caption{Overview of the performance and loss comparison between GEN+GAT, GCN+GAT and GCN+GCN models.}
    \label{fig_supervised_performance}
\end{figure}

\subsection{Results}
The comparison of the performance between 3 models is plotted in Figure \ref{fig_supervised_performance}. We find that GEN+GAT model significantly outperforms the other models GCN+GAT and GCN+GCN. The model GCN+GAT seems to outperform GCN+GCN by few percent, but the performance difference is negligible. \footnote{For GCN+GAT and GCN+GCN models, we have checked that increasing weight dimensions and longer training phases did not lead to better performance.} 

We have also tried other models such as GEN+GEN, GEN+GCN and GAT+GAT to figure out which convolution operators has an important role. The model GEN+GCN shows similar performance with GEN+GAT, but slightly underperforms, and the performance of GAT+GAT is somewhere between that of GEN+GAT and GCN+GAT. This means that GEN plays a significant role to evaluate equivalence or inequivalence for a pair of plumbing graphs. However, we found that GEN+GEN does not perform as good as GEN+GAT or GEN+GCN.

 We used the following datasets to test our models:
 \begin{itemize}
     \item Test set 1\\
     It contains 5,000 pairs of equivalent plumbing graphs generated by \textsc{EquivPair} with $N_{\max}=40$ and 5,000 pairs of inequivalent plumbings generated by \textsc{InequivPair} with $N_{\max}=40$.
     \item Test set 2\\
     This dataset is similar to Test set 1, but with $N_{\max}=60$.
     \item Test set 3\\ 
     This set is also similar to Test set 1, but with $N_{\max}=80$.
     \item Test set 4\\
     It contains 64 pairs of plumbings generated in a manual way such that, for each pair, the determinants of adjacency matrices (with weights on the diagonal) of two plumbings are the same. We use this Test set in order to check that graph embeddings from the models are not just functions of the determinant of the adjacency matrix of a plumbing. All types of Neumann moves have the property that it preserves the determinant of the adjacency matrix of the plumbing, which is the order of the first homology group of the corresponding 3-manifold. We wish graph embeddings to not depend on the determinant only, but be more sophisticated (approximate) invariants of plumbed 3-manifolds.
 \end{itemize}

\begin{figure}[t]
	\centering
	\includegraphics[scale=0.6]{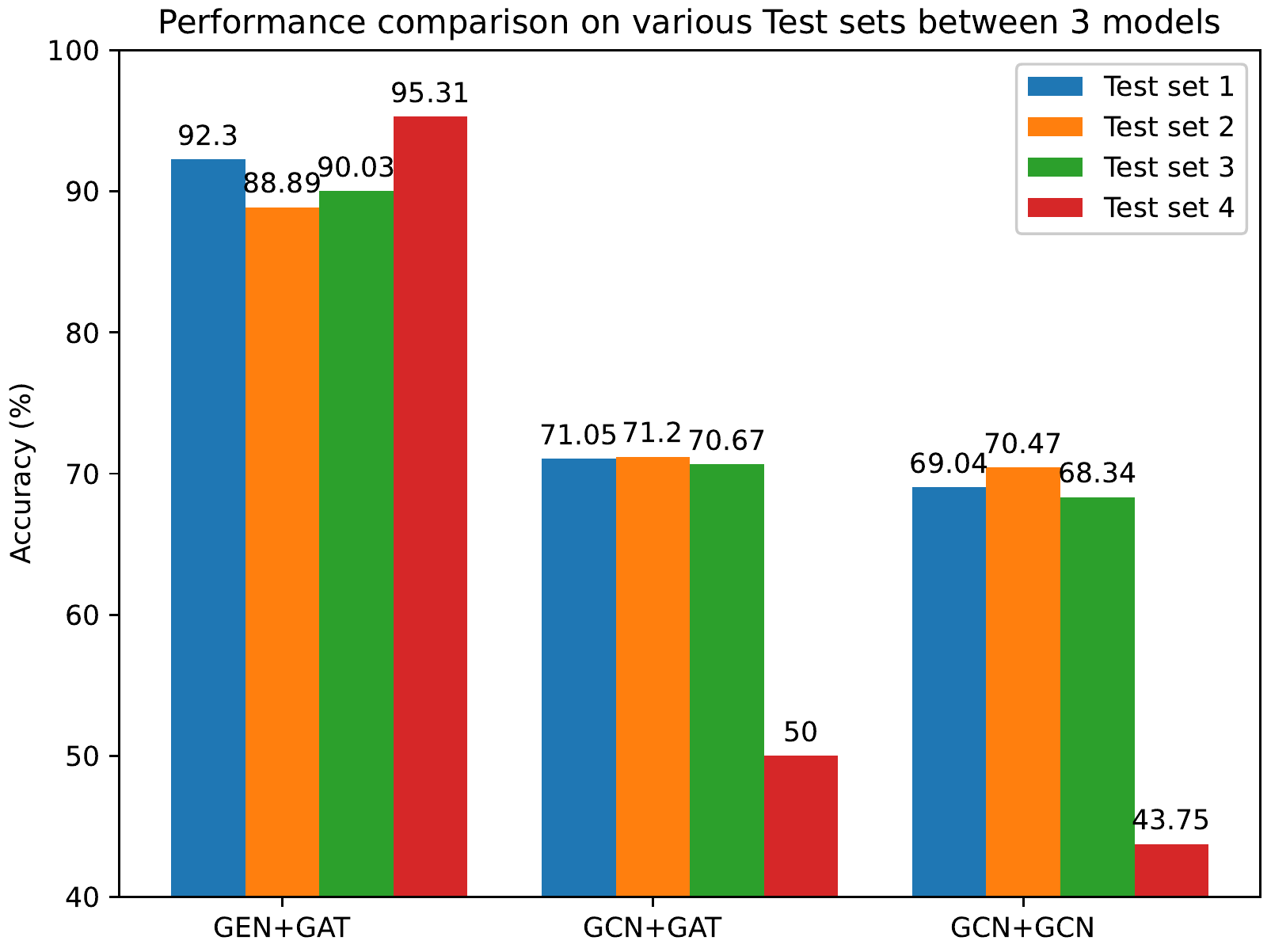}
	\caption{Performance comparison between 3 models, GEN+GAT, GCN+GAT and GCN+GCN on various Test sets. The error bars are not displayed in the figure since the standard errors on Test set 1, 2, and 3 are too small (smaller than $0.7$) to notice. The standard errors on Test set 4 are about 2.64, 6.25, and 6.20 for GEN+GAT, GCN+GAT and GCN+GCN models, respectively.}
	\label{fig_test_result}
\end{figure}

 The results are depicted in Figure \ref{fig_test_result} and they enlighten us with the following two points. The first point is that the accuracy for Test set 2 and Test set 3 is almost the same level as Test set 1 even though Test set 2 and 3 contain plumbing pairs with larger $N_{\max}$ than Test set 1.  It is perhaps surprising that such somewhat counter-intuitive property holds even for GCN+GAT and GCN+GCN models, which show less training accuracy than GEN+GAT. The second point is that GEN+GAT model still outperforms the others for Test set 4 and it can distinguish correctly even inequivalent pairs with the same determinants. Since GEN is designed for graph similarity learning and to have a good generalization, we can see that the model GEN+GAT outperforms significantly the others GCN+GCN and GCN+GAT, designed for general classification problems (with a relatively small number of classes), on various Test sets.

\section{Reinforcement Learning}
\label{section_reinforcement}
In this section, we consider reinforcement learning of a neural network that allows, for a given pair of plumbings, not only to recognize whether they are equivalent or not, but also to find out their simplest representations.

\subsection{The environment}
\subsubsection{State space}
In our RL environment, the plumbing graph defines the state and the state space is infinity. In order to handle the start state and terminal stats in an easy way, we set the start state for an episode is set to be a plumbing generated by \textsc{RandomPlumbing}, Algorithm \ref{alg_randomplumbing}, with number of nodes equal to 10, then applying Neumann moves $N=15$ times.

Between two equivalent plumbings, we define a relation as follows: for two equivalent states $s_1$ and $s_2$, one state is said to be \textit{simpler} than the other if 
\[
    f(s_1) < f(s_2),
\]
where $f(s)$ for a state $s$ is defined by
\begin{equation}
    \label{eqn_simple_definition}
    f(s):= 5|V(s)|+\sum_{v\in V(s)} |w(v)|.
\end{equation}
It is easy to check that this relation is well-defined in a set of all equivalent plumbings. 

\begin{figure}[t]
	\centering
	\includegraphics[scale=0.4]{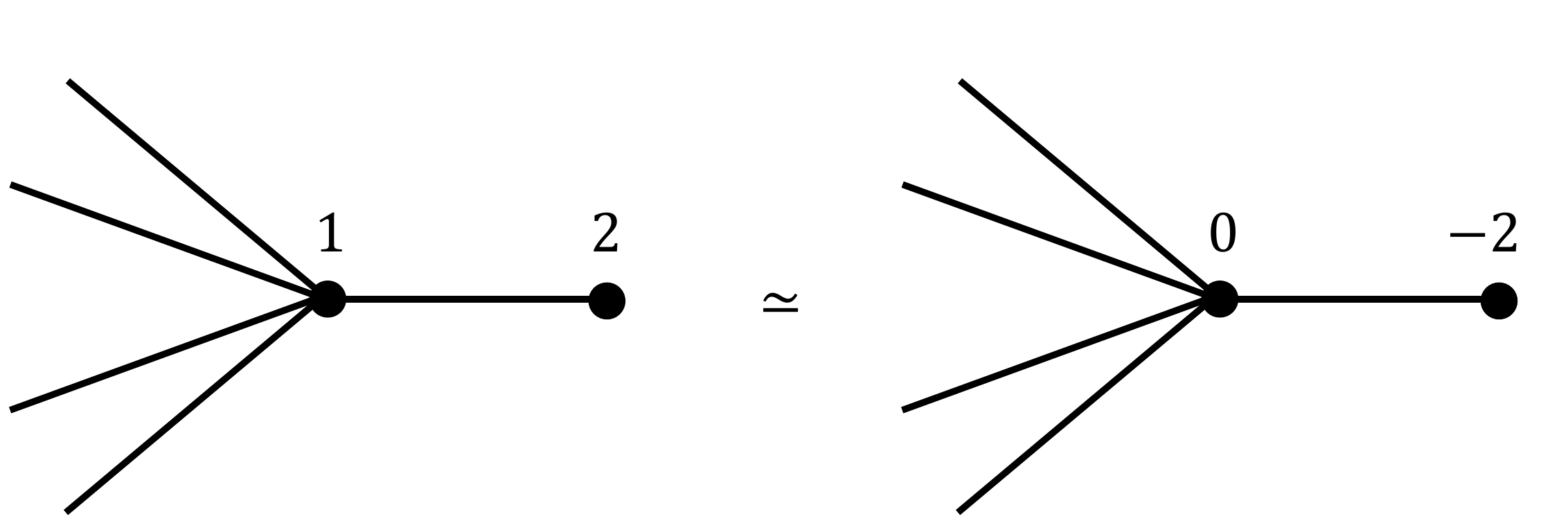}
	\caption{Two plumbings are equivalent and they have same number of nodes. The right-hand side plumbing is simpler than the left one in the sense of \eqref{eqn_simple_definition}.}
	\label{fig_simple_relation}
\end{figure}

One might think that number of nodes in a state is enough to decide which state is simpler. The reason why we add the sum of the absolute values of the weights of nodes is to make the simplest state generically unique\footnote{There still could be specific examples with different plumbings in the same equivalence class that minimize $f(s)$. However, as the results below suggest, such cases are statistically insignificant.}. For example, two plumbings depicted in Figure \ref{fig_simple_relation} have same number of nodes and it is easy to check that they are equivalent by applying 2 Neumann moves. In this example, we say that the plumbing on the right-hand side is simpler than the other from \eqref{eqn_simple_definition}.

By using this comparison relation, we set the terminal state to be a state equal to or simper than the initial state in the episode. We also terminate each episode after taking 15 time steps.

\subsubsection{Action space}
An action for the agent in a state is defined to be a Neumann move applied to one of the nodes. There are 8 possible Neumann moves: 5 blow-up moves and 3 blow-down moves. However, blow-down moves are not always available for all nodes and this could raise a problem that there might be too many of such  \textit{illegal} actions. Therefore, we incorporate 3 blow-down moves into one such that it takes an available blow-down if the corresponding node satisfies one of three following conditions:
\begin{itemize}
    \item the degree of the node is 2 and its weight is equal to $\pm 1$,
    \item the degree of the node is 1 and its weight is equal to $\pm 1$,
    \item the degree of the node is 1 and its weight is equal to 0.
\end{itemize}

Then, for a given state, the total number of possible actions is equal to 6 (5 blow-up moves and 1 blow-down move) times number of nodes in the state. If the agent takes an illegal action, then the next state remains the same state as the current state and the agent will be punished with a negative reward, on which we will elaborate soon.

\subsubsection{Rewards}
Since the goal for the RL agent is to find out the simplest representation for an initial state, it is natural to use $-f(s')$ as a reward (or punishment  $+f(s')$) for taking an action in the current state $s$, where $s'$ denotes the next state obtained by taking an action to the current state $s$. Since all the rewards are negative and simpler state is less punished, it helps the agent not only make the current representation as simple as possible, but also do this job as fast as possible. It is also important to note that some states must get a new blow-up node in order to be simplified, which means the agent has to sacrifice the immediate reward at some time steps to maximize the total return. As we have seen previously, there are some illegal actions in the action space for each state. The reward for such illegal actions is set to be equal to $-2f(s')$ for the next state $s'$, which remains the same as the current state $s$ as we have discussed above.

We set the discount factor as $\gamma = 0.99$, very close to $1$.

\subsection{The deep RL algorithm}
We remind that the RL task is to obtain the simplest representation from a given initial state by using Neumann moves. To accomplish this task, we used Asynchronous Advantage Actor-Critic (A3C) \cite{mnih2016asynchronous} as an RL algorithm, which is the asynchronous version of Actor-Critic (AC) \cite{NIPS1999_6449f44a}, with feedforward GNNs. A3C executes multiple local AC agents asynchronously in parallel to decorrelate the local agent's data into a more stationary process. It also provides practical benefits of being able to use only multi-core CPU, not having to rely on specialized hardware such as GPUs.

The Actor network defines the policy function $\pi(a|s)$, whose output shows the probability of taking action $a$ in state $s$, while the Critic network is to approximate the value function $V^{\pi}(s)$, which represents the expected return from state $s$. Since the inputs of the Actor and Critic are plumbing graphs, in the context of GNNs, the Actor network can be thought as the GNNs for node-level action-selection problem and the Critic is for graph-level estimation problem. The architecture of the Actor is designed by using two graph convolutional layers GCN+GCN and one single-layer feedforward neural network. The Critic has a similar structure, but it has an extra aggregation layer, for which we used a simple mean function. We have also tried GEN+GAT and GEN+GCN for the convolutional layers in the Actor and Critic networks. They seemed to perform well, but it takes a bit longer time for training than GCN+GCN. Since the results with GCN+GCN  were already pretty good, we ended up using GCN+GCN.

We trained the agents for $8\times 10^4$ episodes using 8 CPU cores and no GPU, which takes around 8 hours. We used Adam optimizer with learning rate $5\times 10^{-4}$. For a comparison, we have also implemented Deep Q-Network (DQN) \cite{mnih2013playing} with feedforward GNNs GCN+GCN with the same settings as those for A3C.

\subsection{Results}
Our RL agents can be used to find the simplest representative in the equivalence class of a given plumbing graph. Furthermore, it also can be used to check whether a pair of plumbing graphs represents the same 3-manifold or not. For the latter purpose, we run the RL agents on a pair of plumbings to get the simplest representations for two plumbings, then we compare those to decide whether two equivalent plumbings are isomorphic or not. This process provides us with another advantage that, given two equivalent plumbings, we can get a sequence of Neumann moves that change one plumbing into the other, even though such sequence of Neumann moves is not necessarily the optimal one between two plumbings. From this perspective, we are going to check the performance of the RL agents by running them on pairs of plumbings that represent the same 3-manifolds.  

For the initial inputs of the agents, we generate 10,000 random pairs of plumbings by EQUIVPAIR, Algorithm \ref{alg_equivpair}, but with a fixed number of Neumann moves $N\in\{20, 40, 60, 80, 100\}$. At each time step, the agents choose a Neumann move and apply it to each plumbing in a pair, then we get another pair of plumbings as the next input for the agents. After taking each action, we compare two plumbings and check if they are isomorphic. If yes, we consider it as the success of finding out a sequence of Neumann moves connecting two plumbings in the initial pair. Otherwise, we move on to the next step and we repeat the process until the number of time steps exceeds $5N$. We define the accuracy of the performance as the ratio the number of successes divided by the number of total episodes. An example of a pair of equivalent graphs with the successful result by the A3C trained agent is shown in Figure \ref{figure_two_large_graphs}.

\begin{figure}
     \centering
     \begin{subfigure}[t]{0.49\textwidth}
         \centering
         \includegraphics[width=\textwidth]{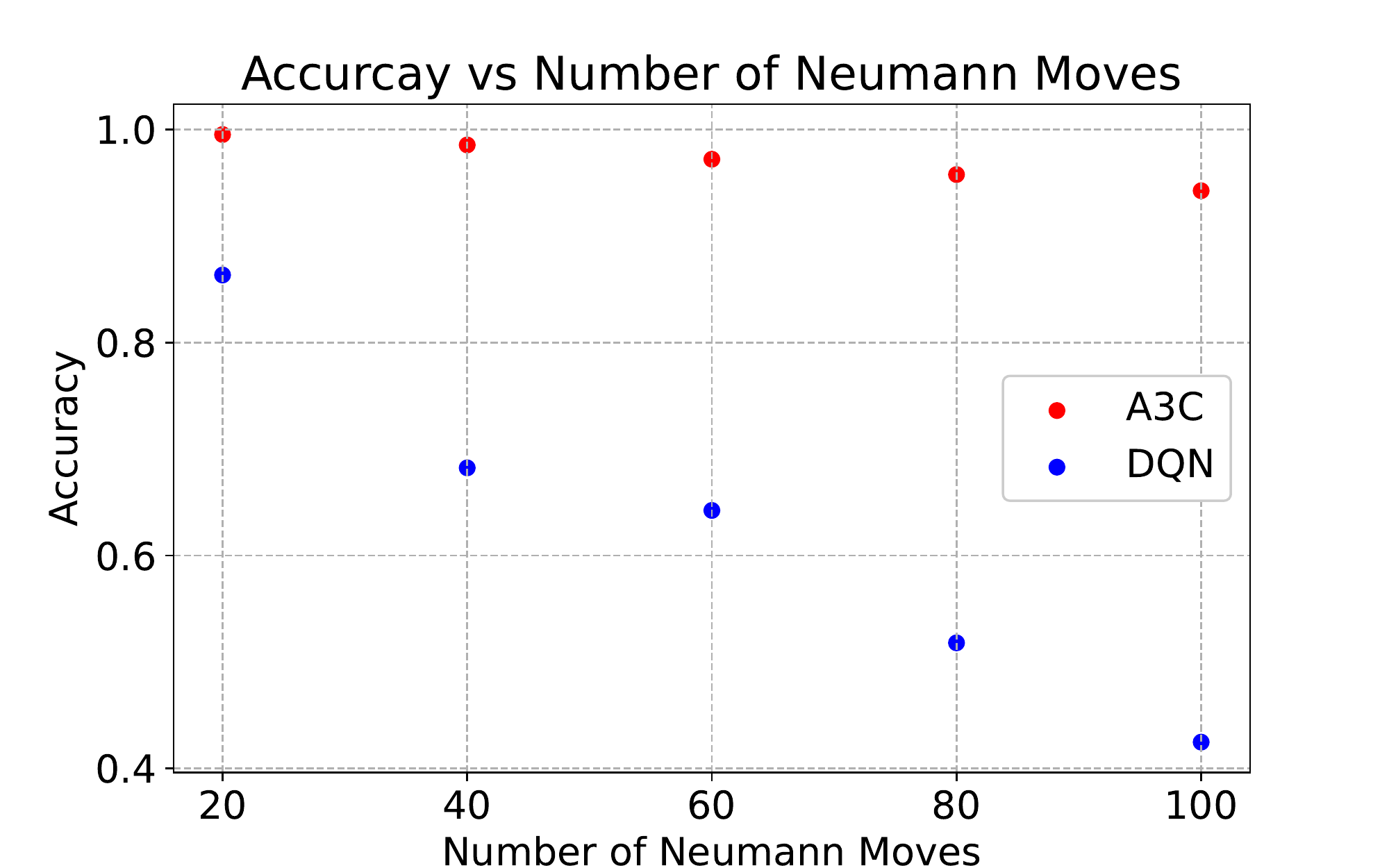}
     \end{subfigure}
     \hfill
     \begin{subfigure}[t]{0.49\textwidth}
         \centering
         \includegraphics[width=\textwidth]{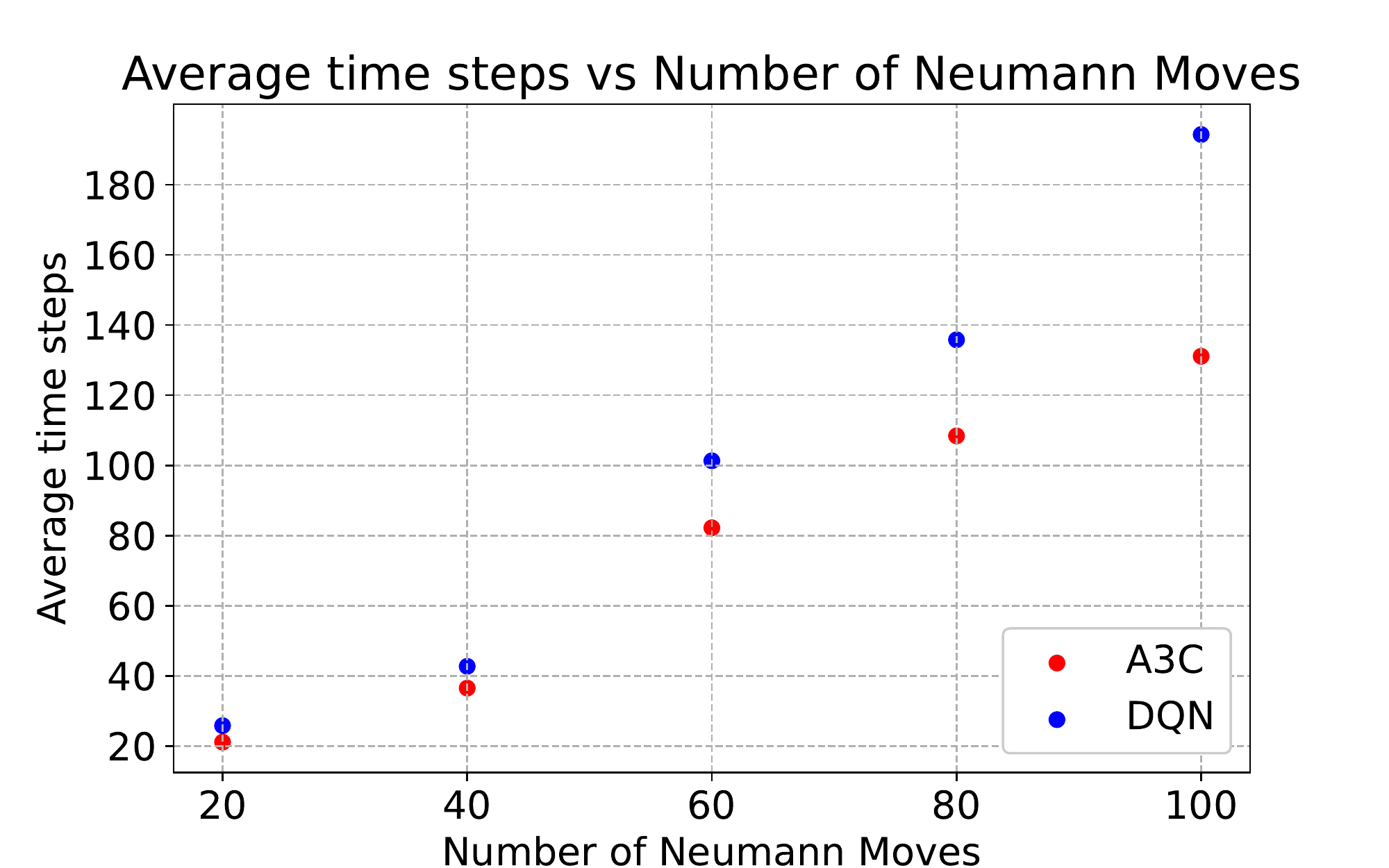}
    \end{subfigure}
    \caption{Performance comparison between A3C and DQN algorithms.}
    \label{fig_rl_test_results}
\end{figure}

The results of the RL agents is presented in Figure \ref{fig_rl_test_results}. The plot on the left in Figure \ref{fig_rl_test_results} shows the accuracy comparison between A3C and DQN. The accuracy for A3C tends to slightly decrease as $N$ gets larger, but it's around 93\% for all pairs of plumbings. However, the accuracy for DQN drops significantly from around 86\% to 42\% when $N$ increases from $N=20$ to $N=100$.

On the right in Figure \ref{fig_rl_test_results}, we show the average number of actions that the agent takes until obtaining a pair of exactly same two plumbings from an initial pair of equivalent plumbings. For A3C agent, the average numbers of actions do not exceed around 1.35 times $N$, which means the trained A3C agent has a good efficiency to make a plumbing simpler. The DQN agent needs similar number of actions to the A3C for $N=20$ and $N=40$. However, it takes almost twice as many number of actions as A3C for larger $N$.

\begin{figure}[t]
	\centering
	\includegraphics[scale=0.7]{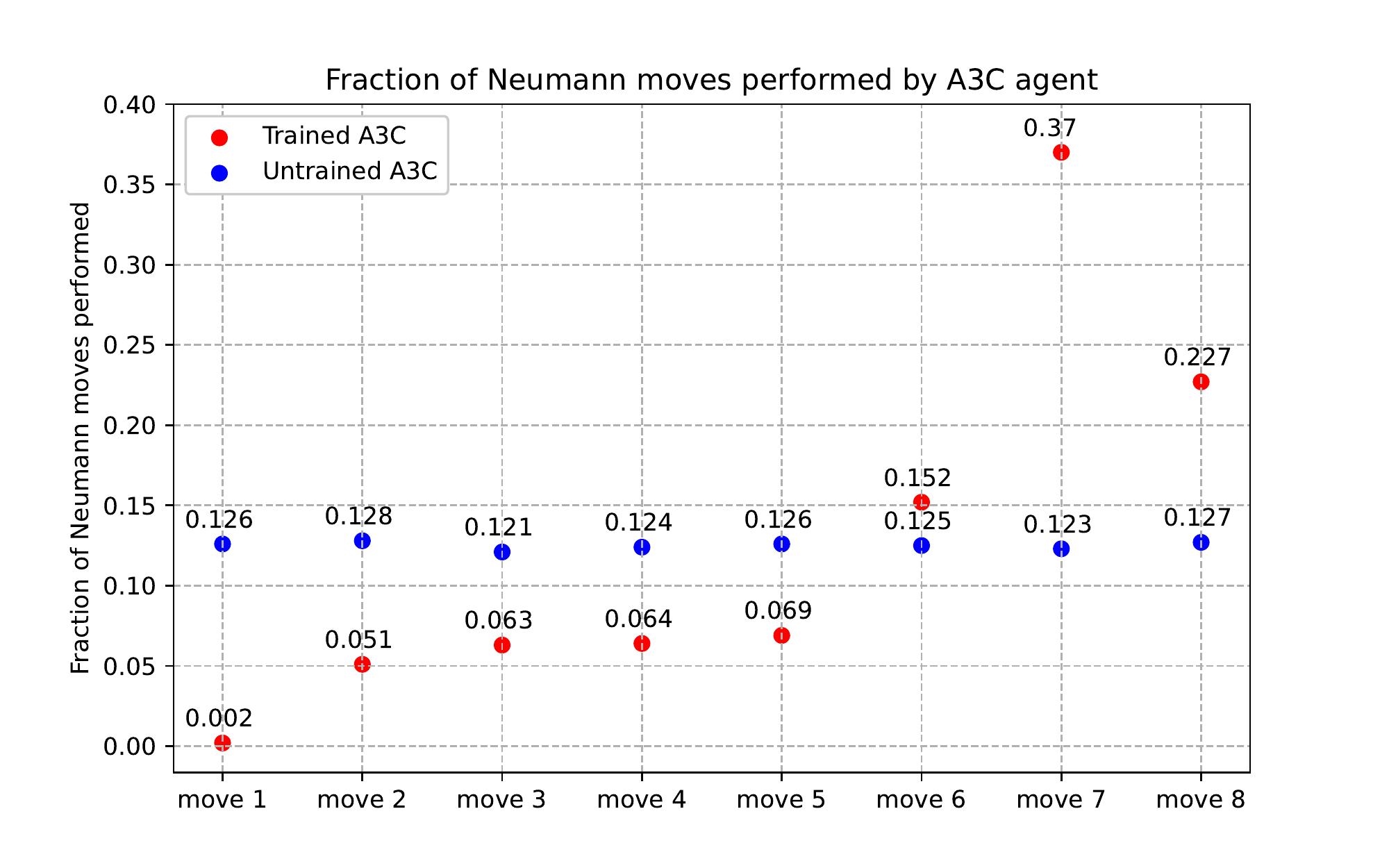}
	\caption{Comparison of the number of Neumann moves taken by a trained A3C agent and an untrained A3C agent to simplify plumbing. The values shown are the total number of Neumann moves of a given type divided by the total number of actions performed, aggregated over multiple examples.}
	\label{fig_statistics_moves}
\end{figure}

We have also studied the distribution of Neumann moves (or actions) that the A3C agent performs before and after training to simplify plumbings generated with $N=100$. In Figure \ref{fig_statistics_moves}, we plot the number of each Neumann move taken by the agent divided by $N$. In the plot, moves 1--5 denote blow-up moves and moves 6--8 denote blow-down moves.

It is natural to observe that all blue dots in Figure \ref{fig_statistics_moves} lay on the line $y=0.125$, because the untrained agent takes each action equally often from a uniform distribution. On the other hand, red dots for trained agent show that the agent takes blow-down moves (moves 6--8) with a probability of around 75\% and takes blow-up moves (moves 1--5) with the remaining probability. This makes sense from the fact that blow-down moves can actually make the plumbing simpler and get a less punishment than blow-up moves. Especially, we see that the move 7, blow-down move of type (b), is the most frequent action and the move 1, blow-up move of type (a), is the least frequent action. This is explained by the fact that the move 1 is not helpful for the agent to get a simpler plumbing.

\begin{figure}[t]
	\centering
	\includegraphics[scale=0.35]{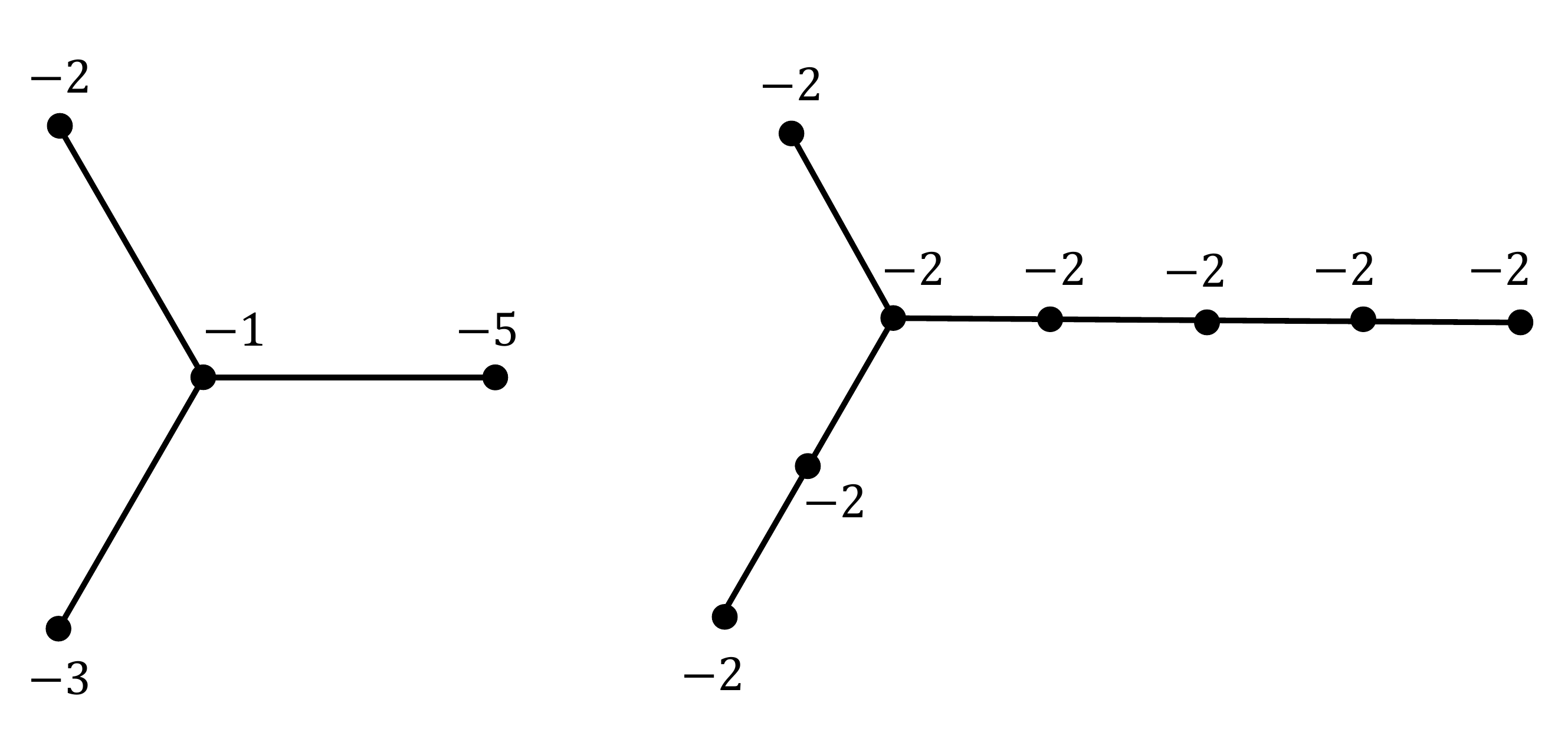}
	\caption{Two equivalent representations of a plumbed 3-manifold $\overline{\Sigma(2,3,5)}$.}
	\label{fig_e8_plumbings}
\end{figure}

Before we jump into the conclusion, it is interesting to check whether or not the trained A3C agent is indeed maximizing the total return instead of immediate rewards by a simple example depicted in \ref{fig_e8_plumbings}. The left plumbing in Figure \ref{fig_e8_plumbings} is a standard representation that realizes a 3-manifold known as a Brieskorn 3-sphere $\overline{\Sigma(2,3,5)}$, while the plumbing on the right represents a homeomorphic 3-manifold which can also be considered as the boundary of the $E_8$ manifold. As one can see immediately, the plumbing on the right in Figure \ref{fig_e8_plumbings} does not have nodes available for blow-down moves. Therefore, in order to get the left plumbing from the right one, the RL agent should take appropriate blow-up moves first, then taking available blow-down moves. This is why we take this example for the test. We notice that 6 actions are needed to turn one plumbing into the other in an optimal way. 

The trained A3C agent successfully simplify the $E_8$ plumbing to the plumbing $\overline{\Sigma(2,3,5)}$ by taking 16 actions, while the trained DQN does not find a solution until the number of actions exceeds 50. This test ensures that the A3C agent indeed pursues not short-term rewards, but its maximal long-term return.

\section{Conclusion and Future Work}
\label{section_conclusion}
\subsection{Conclusion}
In this paper we have examined the GNN approach to the problems in 3-dimensional topology, which ask whether two given plumbing graphs represent a same 3-manifold or not, and whether or not it is possible to find out the sequence of Neumann moves that connects two plumbings if they are equivalent.

In Section \ref{section_supervised}, we used supervised learning to solve the binary classification of whether or not a pair of plumbings is equivalent. We built 3 models by combining graph convolution operators GEN, GCN and GAT, together with a certain graph aggregation module and an MLP as a classifier. We found that GEN+GAT model outperformed GCN+GCN and GCN+GAT models on randomly generated training datasets with maximal number $ N_{\max}=40 $ of applied Neumann moves. GEN+GAT achieved about 95\% accuracy while accuracy for the others is below 80\%. We also tested those 3 models on randomly generated testsets with larger $N_{\max}=60$ and $N_{\max}=80$. Even though those models were trained by a training sets with $N_{\max}=40$, it is an interesting point that, on such testsets, they still performed on a similar level to their training performance.

In Section \ref{section_reinforcement}, we utilized reinforcement learning to find out the sequence of Neumann moves that relates to a given pair of equivalent plumbings. We trained the agent such that it could find the simplest representation of a plumbing by using Neumann moves as its actions. We define the simplicity as a certain linear combination of number of nodes and sum of the absolute value of node features. We ran the trained agent on each of two equivalent plumbings until it arrived at two isomorphic plumbings. In this way, we can construct a sequence of Neumann moves connecting two equivalent plumbings. Using A3C algorithm, we see that the agent can find a sequence of Neumann moves in over 90\% of randomly generated equivalent plumbing pairs even with  $N_{max}=100$. This outperforms the DQN agent by a factor of around 1.5 when $N=60$, and by more than a factor of 2 when $N=100$.

\subsection{Future work}
In this paper we have used Geometric Deep Learning, GNN in particular, in the problem of classification of 3-manifolds up to homeomorphisms. We restricted to a special simple class of 3-manifolds corresponding to tree plumbing graphs. We hope to apply similar neural network models for more general 3-manifolds and also 4-manifolds in the future. One direct generalization would be considering 3-manifolds corresponding to general plumbing graphs described in \cite{neumann1981calculus}, possibly disconnected, with loops, and with non-trivial genera assigned to the vertices\footnote{In such a more general setting, a vertex with weight $w$ and genus $g$ corresponds to a circle fibration of Euler class $w$ over a closed oriented surface of genus $g$. The case considered in this paper is recovered when $g=0$ for each vertex.}. This, in particular, would involve considering extra features associated to the vertices and also to the edges of  graphs, as well as additional set of moves relating equivalent graphs. A more interesting generalization would be considering general Kirby diagrams for 3-manifolds. A Kirby diagram of a 3-manifold is a planar diagram of a link with an integer framing number assigned to each link component. The 3-manifold corresponding to the diagram is then obtained by performing Dehn surgery on this framed link. Two diagrams produce homeomorphic 3-manifolds if and only if they can be related by a sequence of Reidemeister moves (that do not change the isotopy class of the link) together with the so-called Kirby, or equivalently, Fenn-Rourke moves that do change the link but not the resulting 3-manifold (up to homeomorphism). Such a diagram can be understood as a 4-regular plane graph with additional data specifying the types of crossings in the link diagram and the framings of the link components. Alternatively, one can consider Tait graph associated to a checkboard coloring of the link diagram. For practical purposes, this presentation most likely will be more efficient. The Reidemeister, as well as Kirby/Fenn-Rourke moves then can be understood again as certain local operations on graphs associated with Kirby diagrams. The main new challenge would be incorporating the structure of the planar embedding of the graph in GNN. This can be done, for example, by specifying the cyclic order of edges at each vertex, or cyclic order of edges for each face of the plane graph. This additional structure should  be taken into account in the layers of the network. This is not considered in most standard GNN architectures. A further step would be the problem of recognizing whether a pair of Kirby diagrams for 4-manifolds produces a diffeomorphic pair. Such Kirby diagrams are again framed link diagrams that also contain special ``dotted'' link components. There is a corresponding set of local Kirby moves that relate diagrams realizing diffeomorphic 4-manifold. For a comprehensive reference about the Kirby diagrams of 3- and 4-manifold we refer to \cite{gompf19994}.

\section*{Acknowledgements} We would like to thank Sergei Gukov for the useful comments and suggestions on the draft of the paper. We would also like
to thank the anonymous referees who provided insightful and detailed comments and suggestions on a earlier version of the paper.
\appendix
\section{Algorithms}
\label{appendix_algorithms}
In this section, we provide details of the algorithms which have been used to generate datasets for training and testing both SL and RL models in Section \ref{section_supervised} and Section \ref{section_reinforcement}.

\begin{itemize}
    \item \textsc{RandomPlumbing}\\
    This algorithm generates a random plumbing tree by creating a random array for node features and building an adjacency matrix. It starts to choose a random integer as a number of nodes between 1 and 25. In general, there are $N^{N-2} $ different plumbing trees with $N$ nodes if we don't consider node features. Therefore, the upper limit 25 is large enough to generate around $10^6$ random plumbing tress with statistically insignificant overlapping plumbings. The array of node feature is also created by randomly choosing an integer in the interval $(-20, 20)$ for each node. Then we define the adjacency matrix for the plumbing tree, and the algorithm returns a pair of node feature array and adjacency matrix as data for the output plumbing. Note that all random process is done by using a uniform distribution.
    
    \item \textsc{RandomNeumannMove}\\
    The role of this algorithm is to apply a randomly chosen Neumann move to a random node of the input plumbing, then returns the resulting plumbing. A random Neumann move is characterized by 3 variables, i.e., $type$, $updown$, and $sign$. Here $type\in \{1,2,3\}$ denotes 3 types of Neumann moves depicted in Figure \ref{fig_nmoves}, $updown\in \{1, -1\}$ points out blow-up ($updown = 1$) or blow-down ($updown=-1$), and $sign\in \{1, -1\}$ denotes the sign of the new vertex for blow-up Neumann moves of type (b) and (c). Notice that other moves does not require $sign$.

    The algorithm first takes a random node of the input and fixes a random tuple $(type, updown, sign)$ from a uniform distribution. Then it builds new node feature array and adjacency matrix for the plumbing obtained by applying the Neumann move to the chosen node. If the Neumann move determined by a tuple $(type, updown, sign)$ is an illegal move, the output plumbing is the same as the input. The algorithm also returns another variable $done\in \{\textsc{True, False}\}$, which makes it possible to notice whether the Neumann move to be applied is legal ($done=\textsc{True}$) or illegal ($done=\textsc{False}$).  This variable $done$ will be used to decide the rewards of actions in Section \ref{section_reinforcement}.

    \item \textsc{EquivPair} and \textsc{InequivPair}\\
    These are used to generate an equivalent plumbing pair (\textsc{EquivPair}) or an inequivalent plumbing pair (\textsc{InequivPair}). At the first step, \textsc{EquivPair} generates an initial pair of isomorphic plumbings, while \textsc{InequivPair} generates two inequivalent plumbings, by using \textsc{RandomPlumbing}. Then they have the same process, in which they apply Neumann moves iteratively up to $N_{\max}$ times to each plumbing in the initial pair. Then they return the resulting pair as well as a variable, named $label$, which will be used for classification problem in Section \ref{section_supervised}. Notice that $label = 1$ for \textsc{EquivPair} and $label = -1$ for \textsc{InequivPair}. 
    
    \item \textsc{TweakPair}\\
    This algorithm generates an inequivalent pair of plumbings, but with the same graph structure. One plumbing is generated by \textsc{RandomPlumbing}, and the other is obtained by tweaking a copy of the first plumbing, i.e., by making a small change to a feature of a randomly chosen node. These two plumbings form an initial pair. Since the adjacency matrices of two plumbings are same, they have the same graph structure. However, due to the small change, two plumbings are inequivalent. Then the algorithm has the same structure as in \textsc{EquivPair} and \textsc{InequivPair} to apply random Neumann moves iteratively to each plumbing in the initial pair.
\end{itemize}

\begin{algorithm}
    \caption{\textsc{RandomPlumbing}}
    \label{alg_randomplumbing}
    \begin{algorithmic}
    \State $n \gets$ random integer between 1 and 25 \Comment{number of nodes}
    \State $\mathbf{x} \gets$ array of $n$ random integers between $-20$ and 20 \Comment{node features}
    \State $\mathbf{a} \gets$ $n\times n$ matrix of zeros \Comment{initialize the adjacency matrix}
    \For{\textnormal{$i = 2$ to $n$}}\Comment{construct the adjacency matrix}
        \State {$j \gets$ random integer between 1 and $i-1$}
        \State {$\mathbf{a}_{i,j}, \mathbf{a}_{j,i} \gets 1$}
    \EndFor
    \State $G\gets$ $(\mathbf{x}, \mathbf{a})$ \Comment{$G$ defines the plumbing}
    \State \textbf{return} $G$
    \end{algorithmic}
\end{algorithm}    
    
\begin{algorithm}
\caption{\textsc{RandomNeumannMove}}
\label{alg_RandomNeumannMove}
\begin{algorithmic}
    \State \textbf{Require:} a plumbing $G$
    \State $v\gets$ a random node of $G$
    \State $type\gets$ a random choice in $\{1,2,3\}$
    \State $updown\gets$ a random choice in $\{1, -1\}$
    \If{$updown=1$}\Comment{blow-up move}
        \If {$type=1$}
            \State $G'\gets$ a plumbing applied a blow-up move of type (a) to the node $v$
        \Else
            \State $sign\gets$ a random choice in $\{1, -1\}$
            \State $G'\gets$ a plumbing applied a blow-up move determined by $(type, sign)$
        \EndIf
        \State $done\gets$ \textsc{True}
    \Else\Comment{blow-down move}
        \If{$v$ can be removed by a blow-down move}
            \State $G'\gets$ a plumbing applied a blow-down move to the node $v$
            \State $done\gets$ \textsc{True}
        \Else
            \State $G'\gets G$ \Comment{returns the input plumbing for a forbidden move}
            \State $done\gets$ \textsc{False}
        \EndIf
    \EndIf
    \State \textbf{return} $(done, G')$           
\end{algorithmic}
\end{algorithm}

\begin{algorithm}[t]
    \caption{EQUIVPAIR}\label{alg_equivpair}
    \begin{algorithmic}
        \State \textbf{Require:} $N_{\max}\in \mathbb{Z}^{+}$
        \State $G\gets$ a plumbing by \textsc{RandomPlumbing}
        \State $G_1\gets G$
        \State $n_1\gets$ a random integer between 1 and $N_{\max}$
        \For {$i=1$ to $n_1$}\Comment{Apply Neumann moves $n_1$ times}
            \State $G_1\gets$ \textsc{RandomNeumannMove}$(G_1)$
        \EndFor
        \State $G_2\gets G$
        \State $n_2\gets$ a random integer between 1 and $N_{\max}$
        \For {$j=1$ to $n_2$}\Comment{Apply Neumann moves $n_2$ times}
            \State $G_2\gets$ \textsc{RandomNeumannMove}$(G_2)$
        \EndFor
        \State $label\gets 1$
        \State \textbf{return} $G_1, G_2, label$
    \end{algorithmic}
\end{algorithm}

\begin{algorithm}[t]
    \caption{INEQUIVPAIR}\label{alg_inequivpair}
    \begin{algorithmic}
    \State \textbf{Require:} $N_{\max}\in \mathbb{Z}^{+}$
        \State $G_1\gets$ a plumbing by \textsc{RandomPlumbing}
        \State $G_2\gets$ another plumbing by \textsc{RandomPlumbing}
        \State $n_1\gets$ a random integer between 1 and $N_{\max}$
        \For {$i=1$ to $n_1$}\Comment{Apply Neumann moves $n_1$ times}
            \State $G_1\gets$ \textsc{RandomNeumannMove}$(G_1)$
        \EndFor
        \State $n_2\gets$ a random integer between 1 and $N_{\max}$
        \For {$j=1$ to $n_2$}\Comment{Apply Neumann moves $n_2$ times}
            \State $G_2\gets$ \textsc{RandomNeumannMove}$(G_2)$
        \EndFor
        \State $label\gets -1$
        \State \textbf{return} $G_1, G_2, label$
    \end{algorithmic}
\end{algorithm}

\begin{algorithm}[t]
    \caption{\textsc{TweakPair}}\label{alg_tweakpair}
    \begin{algorithmic}
        \State \textbf{Require:} $N_{\max}\in \mathbb{Z}^{+}$
        \State $G_1\gets$ a plumbing by \textsc{RandomPlumbing}
        \State $G_2\gets G_1$ 
        \State $v\gets$ a random node in $G_2$
        \State $t\gets$ a random integer between -3 and 3, not 0.
        \State $\mathbf{x}\gets$ node feature of $G_2$
        \State $\mathbf{a}\gets$ adjacency matrix of $G_2$
        \State $\mathbf{x}_v\gets \mathbf{x}_v+t$
        \State $G_2 \gets$ a plumbing with $(\mathbf{x}, \mathbf{a})$
        
        \State $n_1\gets$ a random integer between 1 and $N_{\max}$
        \For {$i=1$ to $n_1$}\Comment{Apply Neumann moves $n_1$ times}
            \State $G_1\gets$ \textsc{RandomNeumannMove}$(G_1)$
        \EndFor
        \State $n_2\gets$ a random integer between 1 and $N_{\max}$
        \For {$i=1$ to $n_2$}\Comment{Apply Neumann moves $n_2$ times}
            \State $G_2\gets$ \textsc{RandomNeumannMove}$(G_2)$
        \EndFor
        \State $label\gets -1$
        \State \textbf{return} $G_1, G_2, label$
    \end{algorithmic}
\end{algorithm}    

\begin{figure}[b]
	\centering
	\includegraphics[scale=0.5]{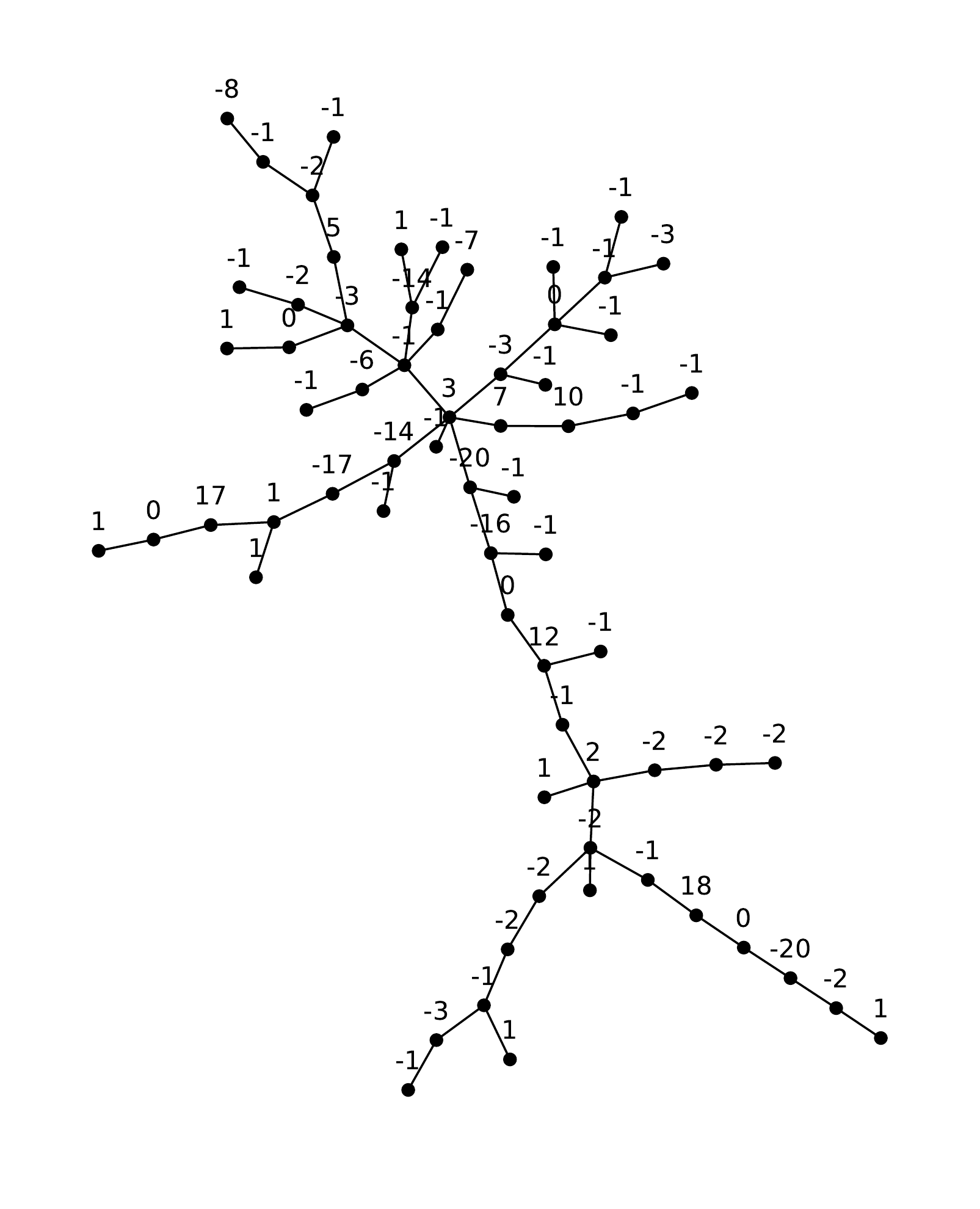}
	\includegraphics[scale=0.5]{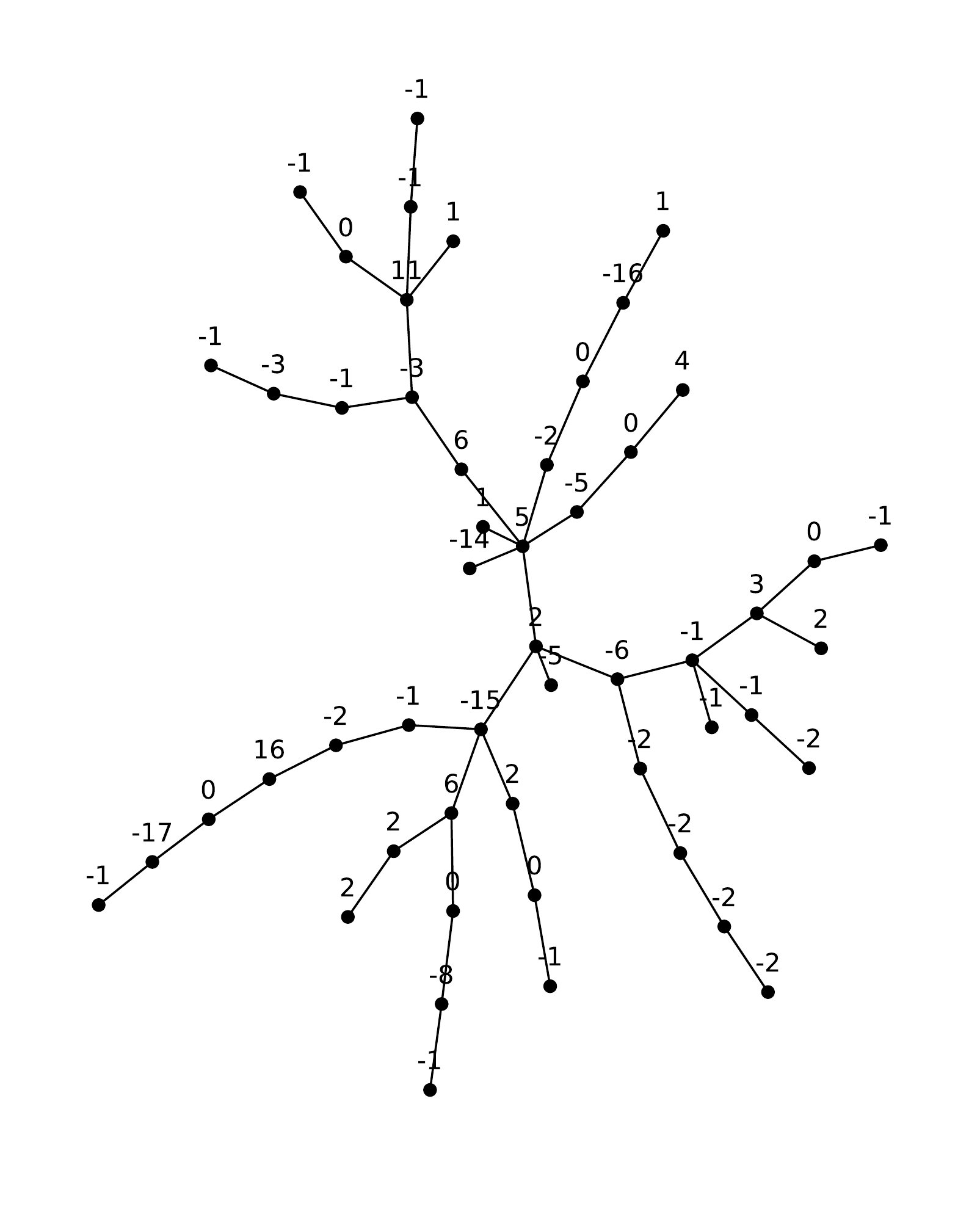}
 \caption{An example of a pair of equivalent plumbing graphs generated by EQUIVPAIR with the number of Neumann moves fixed to $N=40$. The graphs are successfully recognized as equivalent both by the RL agent trained by A3C algorithm considered in Section \ref{section_reinforcement} and the GEN+GAT neural network considered in Section \ref{section_supervised}. }
 \label{figure_two_large_graphs}
\end{figure}
\clearpage
\bibliographystyle{JHEP}
\bibliography{refs}

\providecommand{\href}[2]{#2}\begingroup\raggedright\begin{thebibliography}{10}

\bibitem{bronstein2021geometric}
M.M.~Bronstein, J.~Bruna, T.~Cohen and P.~Veličković, \emph{Geometric deep
  learning: Grids, groups, graphs, geodesics, and gauges},  2021.

\bibitem{cao2020comprehensive}
W.~Cao, Z.~Yan, Z.~He and Z.~He, \emph{A comprehensive survey on geometric deep
  learning}, {\emph{IEEE Access} {\bfseries 8} (2020) 35929}.

\bibitem{hughes2020neural}
M.C.~Hughes, \emph{A neural network approach to predicting and computing knot
  invariants}, {\emph{Journal of Knot Theory and Its Ramifications} {\bfseries
  29} (2020) 2050005}.

\bibitem{Jejjala:2019kio}
V.~Jejjala, A.~Kar and O.~Parrikar, \emph{{Deep Learning the Hyperbolic Volume
  of a Knot}},
  \href{https://doi.org/10.1016/j.physletb.2019.135033}{\emph{Phys. Lett. B}
  {\bfseries 799} (2019) 135033}
  [\href{https://arxiv.org/abs/1902.05547}{{\ttfamily 1902.05547}}].

\bibitem{Gukov:2020qaj}
S.~Gukov, J.~Halverson, F.~Ruehle and P.~Su\l{}kowski, \emph{{Learning to
  Unknot}}, \href{https://doi.org/10.1088/2632-2153/abe91f}{\emph{Mach. Learn.
  Sci. Tech.} {\bfseries 2} (2021) 025035}
  [\href{https://arxiv.org/abs/2010.16263}{{\ttfamily 2010.16263}}].

\bibitem{davies2021advancing}
A.~Davies, P.~Veli{\v{c}}kovi{\'c}, L.~Buesing, S.~Blackwell, D.~Zheng,
  N.~Toma{\v{s}}ev et~al., \emph{Advancing mathematics by guiding human
  intuition with ai}, {\emph{Nature} {\bfseries 600} (2021) 70}.

\bibitem{kauffman2022rectangular}
L.~Kauffman, N.~Russkikh and I.~Taimanov, \emph{Rectangular knot diagrams
  classification with deep learning}, {\emph{Journal of Knot Theory and Its
  Ramifications} {\bfseries 31} (2022) 2250067}.

\bibitem{Craven:2021ckk}
J.~Craven, M.~Hughes, V.~Jejjala and A.~Kar, \emph{{Learning knot invariants
  across dimensions}},
  \href{https://doi.org/10.21468/SciPostPhys.14.2.021}{\emph{SciPost Phys.}
  {\bfseries 14} (2023) 021}
  [\href{https://arxiv.org/abs/2112.00016}{{\ttfamily 2112.00016}}].

\bibitem{vernitski2022reinforcement}
A.~Vernitski, A.~Lisitsa et~al., \emph{Reinforcement learning algorithms for
  the untangling of braids},  in \emph{The International FLAIRS Conference
  Proceedings}, vol.~35, 2022.

\bibitem{khan2021untangling}
A.~Khan, A.~Vernitski and A.~Lisitsa, \emph{Untangling braids with multi-agent
  q-learning},  in \emph{2021 23rd International Symposium on Symbolic and
  Numeric Algorithms for Scientific Computing (SYNASC)}, pp.~135--139, IEEE,
  2021.

\bibitem{lisitsa2023supervised}
A.~Lisitsa, M.~Salles and A.~Vernitski, \emph{Supervised learning for
  untangling braids},  in \emph{15th International Conference on Agents and
  Artificial Intelligence (ICAART 2023)}, 2023.

\bibitem{gukov2023searching}
S.~Gukov, J.~Halverson, C.~Manolescu and F.~Ruehle, \emph{Searching for ribbons
  with machine learning},  2023.

\bibitem{He:2023csq}
Y.-H.~He, E.~Heyes and E.~Hirst, \emph{{Machine Learning in Physics and
  Geometry}},  \href{https://arxiv.org/abs/2303.12626}{{\ttfamily 2303.12626}}.

\bibitem{neumann1981calculus}
W.D.~Neumann, \emph{A calculus for plumbing applied to the topology of complex
  surface singularities and degenerating complex curves}, {\emph{Transactions
  of the American Mathematical Society} {\bfseries 268} (1981) 299}.

\bibitem{kuperberg2019algorithmic}
G.~Kuperberg, \emph{Algorithmic homeomorphism of 3-manifolds as a corollary of
  geometrization}, {\emph{Pacific Journal of Mathematics} {\bfseries 301}
  (2019) 189}.

\bibitem{ZHOU202057}
J.~Zhou, G.~Cui, S.~Hu, Z.~Zhang, C.~Yang, Z.~Liu et~al., \emph{Graph neural
  networks: A review of methods and applications},
  \href{https://doi.org/https://doi.org/10.1016/j.aiopen.2021.01.001}{\emph{AI
  Open} {\bfseries 1} (2020) 57}.

\bibitem{li2019graph}
Y.~Li, C.~Gu, T.~Dullien, O.~Vinyals and P.~Kohli, \emph{Graph matching
  networks for learning the similarity of graph structured objects},  2019.

\bibitem{kipf2017semisupervised}
T.N.~Kipf and M.~Welling, \emph{Semi-supervised classification with graph
  convolutional networks},  2017.

\bibitem{velickovic2018graph}
P.~Veličković, G.~Cucurull, A.~Casanova, A.~Romero, P.~Liò and Y.~Bengio,
  \emph{Graph attention networks},  in \emph{International Conference on
  Learning Representations}, 2018,
  \href{https://openreview.net/forum?id=rJXMpikCZ}{https://openreview.net/forum?id=rJXMpikCZ}.

\bibitem{paszke2017pytorch}
A.~Paszke, S.~Gross, S.~Chintala, G.~Chanan, E.~Yang, Z.~DeVito et~al.,
  \emph{Automatic differentiation in pytorch},  in \emph{NIPS-W}, 2017.

\bibitem{FeyLenssen2019pytorchgeometric}
M.~Fey and J.E.~Lenssen, \emph{Fast graph representation learning with {PyTorch
  Geometric}},  in \emph{ICLR Workshop on Representation Learning on Graphs and
  Manifolds}, 2019,
  \href{https://arxiv.org/abs/1903.02428}{https://arxiv.org/abs/1903.02428}.

\bibitem{li2017gated}
Y.~Li, D.~Tarlow, M.~Brockschmidt and R.~Zemel, \emph{Gated graph sequence
  neural networks},  2017.

\bibitem{mnih2016asynchronous}
V.~Mnih, A.P.~Badia, M.~Mirza, A.~Graves, T.P.~Lillicrap, T.~Harley et~al.,
  \emph{Asynchronous methods for deep reinforcement learning},  2016.

\bibitem{NIPS1999_6449f44a}
V.~Konda and J.~Tsitsiklis, \emph{Actor-critic algorithms},  in \emph{Advances
  in Neural Information Processing Systems}, S.~Solla, T.~Leen and
  K.~M\"{u}ller, eds., vol.~12, MIT Press, 1999,
  \href{https://proceedings.neurips.cc/paper\_files/paper/1999/file/6449f44a102fde848669bdd9eb6b76fa-Paper.pdf}{https://proceedings.neurips.cc/paper\_files/paper/1999/file/6449f44a102fde848669bdd9eb6b76fa-Paper.pdf}.

\bibitem{mnih2013playing}
V.~Mnih, K.~Kavukcuoglu, D.~Silver, A.~Graves, I.~Antonoglou, D.~Wierstra
  et~al., \emph{Playing atari with deep reinforcement learning},  2013.

\bibitem{gompf19994}
R.E.~Gompf and A.~Stipsicz, \emph{4-manifolds and Kirby calculus}, vol.~20 of
  \emph{Graduate Studies in Mathematics}, American Mathematical Soc. (1999).

\end{thebibliography}\endgroup
%----------------------------------------------------------------------------------------

\end{document}